# Iterative polynomial-root-finding procedure with enhanced accuracy

Hashim A. Yamani [a] and Abdulaziz D. Alhaidari [b] [†]

[a] Mabuth-2414, Medina 42362-6959, Saudi Arabia

[b] Saudi Center for Theoretical Physics, P.O. Box 32741, Jeddah 21438, Saudi Arabia

**Abstract**. We devise a simple but remarkably accurate iterative routine for calculating the roots of a polynomial of any degree. We demonstrate that our results have significant improvement in accuracy over those obtained by methods used in popular software packages.

**It turned out that the proposed procedure in this work is a rediscovery of the "Durand-Kerner" method. See, for example,**

E. Durand, *Solutions numériques des équations algébriques. Tome I, Équations du type F(x) = 0, racines d'un polynôme*, Masson, Paris, 1960.

K. Dochev, *Modified Newton's method for simultaneous computation of all the roots of a given algebraic equation* (Bulgarian), Phys. Mat. J. Bulg. Acad. Sci. 5 (1962), 136-139.

I. O. Kerner, *Ein gesamtschrittverfahren zur berechnung der nullstellen von polynomen*, Numerische Mathematik 8 (1966) 290-294.

## 1. Introduction

The search for algebraic formulas for the zeros of higher order polynomial equations stopped at the fourth order (the quartic equation) because in 1820 the Norwegian mathematician Abel showed that no such formula exists for degrees higher than 4 [1]. Thus, numerical routines for finding the roots of polynomials of degrees higher than four were needed and the race to come up with the most efficient routine (in convergence, accuracy, speed and stability) was on and still going. In this work, we develop a highly accurate iterative procedure for calculating the roots of a polynomial of any degree, which we write in monic form as follows

$$P_n(z) = z^n + a_{n-1}z^{n-1} + a_{n-2}z^{n-2} + ... + a_1 z + a_0 = z^n + \sum_{m=0}^{n-1} a_m z^m , \qquad (1)$$

---

[†] Corresponding Author: email address: haidari@sctp.org.sa

where the $n$ coefficients $\{a_m\}_{m=0}^{n-1}$ are complex numbers. The fundamental theorem of algebra dictates that this polynomial will have $n$ complex roots (zeros). If we call these roots $\{z_m\}_{m=0}^{n-1}$, then we can factorize the polynomial as follows

$$P_n(z) = (z-z_0)(z-z_1)...(z-z_{n-1}) = \prod_{m=0}^{n-1}(z-z_m). \qquad (2)$$

Methods, such as the New-Raphson, Brent, Laguerre, and their variants, have been devised to search and find the roots $\{z_m\}_{m=0}^{n-1}$ starting from the known coefficients $\{a_m\}_{m=0}^{n-1}$. These methods invariably adopt a two-element strategy [2]. The first is the use of iteration starting with an initial approximate guess of one of the roots. Convergence to an accurate value of the sought root depends on how good the initial guess is. Methods differ in the details of how to pick the initial guess. The second element is using staged deflation of the polynomial. This occurs by factoring out the term involving the converged root thereby reducing the polynomial by one degree. This process is applied to each resulting polynomial. The final stage is reached when only the last root remains to be found. Then, we can obtain an accurate value of this root by rewriting (2) as

$$z_0 = z - \frac{P_n(z)}{(z-z_1)(z-z_2)...(z-z_{n-1})}, \qquad (3)$$

and evaluating the right-hand side at any arbitrary point $\tilde{z}_0$ in the complex plane (of course, other than the points $\{z_m\}_{m=1}^{n-1}$). However, if these roots are known only approximately as $\{\hat{z}_m\}_{m=1}^{n-1}$ then obviously the accuracy of the resulting $z_0$ depends on those of $\{\hat{z}_m\}_{m=1}^{n-1}$ and our choice of the arbitrary point $\tilde{z}_0$.

In our approach, we also adopt the same two-elements strategy, but as will be shown in the next section, we exploit the power of deflation of the polynomial to its fullest. Instead of the staged deflation of the polynomial every time an accurate root is found, we propose to fully deflate the polynomial using approximate roots simultaneously and from the start of the iteration procedure

## 2. Procedure

The process of finding roots of a polynomial suggests a natural definition of accuracy of the process. We take as measure of the accuracy of any given root $z_k$ to be the smallness of the value $|P_n(z_k)|$. That is, a calculation of this root that gives the value $z_k'$ is more accurate than another that gives $z_k''$ if $|P_n(z_k')| < |P_n(z_k'')|$.

Our procedure starts out with arbitrary initial guess values $\{\hat{z}_m\}_{m=0}^{n-1}$. Then, we make an iterative calculation of $z_0$ using (3) as follows

$$s_{k+1} = s_k - \frac{P_n(s_k)}{(s_k - \hat{z}_1)(s_k - \hat{z}_2)\ldots(s_k - \hat{z}_{n-1})}, \qquad (4)$$

for $k = 0,1,2,\ldots$ and with $s_0$ being the arbitrary point $\hat{z}_0$. Then for a given precision (number of accurate digits) the calculation converges after a finite number of iterations $N$ to a fixed value, which we call $\bar{z}_0 = s_N$. Subsequently, we proceed to improve $\hat{z}_1$ iteratively using the accurately calculated $\bar{z}_0$ and writing (4) for this root as follows

$$s_{k+1} = s_k - \frac{P_n(s_k)}{(s_k - \bar{z}_0)(s_k - \hat{z}_2)(s_k - \hat{z}_3)\ldots(s_k - \hat{z}_{n-1})}, \qquad (5)$$

for $k = 0,1,2,\ldots$ and with $s_0 = \hat{z}_1$. The calculation will then converge to a fixed value, which we call $\bar{z}_1$. We repeat the same for $\hat{z}_2$ using $\bar{z}_0$ and $\bar{z}_1$ by writing

$$s_{k+1} = s_k - \frac{P_n(s_k)}{(s_k - \bar{z}_0)(s_k - \bar{z}_1)(s_k - \hat{z}_3)(s_k - \hat{z}_4)\ldots(s_k - \hat{z}_{n-1})}, \qquad (6)$$

for $k = 0,1,2,\ldots$ and with $s_0 = \hat{z}_2$. And so on, where an improvement of the approximate root $\hat{z}_m$ is obtained from the convergent iteration

$$s_{k+1} = s_k - \frac{P_n(s_k)}{\prod_{j=0}^{m-1}(s_k - \bar{z}_j)\prod_{j=m+1}^{n-1}(s_k - \hat{z}_j)}. \qquad (7)$$

Now, if the desired accuracy of the resulting roots $\{\bar{z}_m\}_{m=0}^{n-1}$ is not reached, as measured by how small the set $\{|P_n(\bar{z}_m)|\}_{m=0}^{n-1}$ is, then we repeat the procedure for another round but now we start with the new improved values $\{\bar{z}_m\}_{m=0}^{n-1}$ rather than the guessed values $\{\hat{z}_m\}_{m=0}^{n-1}$. And so on[‡].

In the following section, we apply the procedure to an example of a polynomial of degree 20 with arbitrarily chosen complex coefficients and compare the accuracy of our results to those obtained by other well-known computing software packages.

## 3. Results and discussion

As illustration, we used the computational software package Mathcad® to write the routine **Khandug** implementing the procedure in the previous section and as shown in Figure 1. We compare the accuracy of the routine to that of the built-in **polyroots** function in Mathcad, which utilizes either the "La Guerre (LG)" method or the "Companion Matrix (CM)" method [3]. We also compare with the results obtained by the three built-in routines **NRoots**, **Nsolve** and **Reduce** in the computational software package Mathematica®. We took $n = 20$ and made an arbitrary choice of the complex coefficients $\{a_m\}_{m=0}^{19}$ shown in the first column of Table I. We have chosen the initial guessed values for the roots $\{\hat{z}_m\}_{m=0}^{19}$ as 20 equally spaced points on the unit circle in the complex plane. The routine executed one iteration step given by Eq. (7) for each root (i.e., $N = 1$) and then the whole procedure was repeated 17 rounds to reach the desired accuracy. Table II shows a comparison of the accuracy of the roots obtained by our iterative routine **Khandug** to those obtained by Mathcad and Mathematica. The first column in Table II is a measure of how good our initial guessed values by giving $\{|P_n(\hat{z}_m)|\}_{m=0}^{19}$. It is clear from the Table that the accuracy of our results is superior to those obtained by the two computational packages. Other tests do also confirm the preeminence of our routine [4]. We also compared our results to those obtained by two FORTRAN routines from the SLATEC Common Mathematical Library [5], **CPZERO** and **CPQE79**, with the same conclusion that **Khandug** is not just more accurate but also faster [6]. It should be mentioned that we could have also reached the same desired accuracy as in Table II by performing 30 iterations for each root ($N = 30$) and then repeating the procedure for three rounds instead of 17. However, this takes a total of 1800 iteration step whereas for Table II it is only 340. We show in Table III the lowest iteration sequences in which the process converges to the desired accuracy for the same initial guessed roots.

---

[‡] Beside the accuracy measure for each individual root, $|P_n(z_m)|$, there are overall (global) accuracy measures for the entire root-finding procedure. For example, we can utilize the well-known relations of the sum and product of the roots to the coefficients $a_{n-1}$ and $a_0$, respectively. Therefore, we can use as global measure of accuracy the smallness of either $\left|a_{n-1} + \sum_{m=0}^{n-1} z_m\right|$ or $\left|(-1)^n a_0 - \prod_{m=0}^{n-1} z_m\right|$. These global accuracy measures could be used as the terminating condition for iterations in the root-finding routine.

In the last row of Table II, it is peculiar to see that Mathcad and Mathematica produced a relatively large number for one of the roots, which is nine to twelve orders of magnitudes larger than the largest of the other numbers in the corresponding column. These large numbers correspond to the root "−0.09074+5.81549i" whose accuracy in our calculation is even better than $10^{-307}$. The same peculiar result is produced by other computational tools [4]. This result is not a precision-related issue since it persists even after increasing precision to the maximum allowed by our computing hardware and software. Of course, if we repeat the calculation with increased precision, the value for this number decreases but its relatively large size to the rest of the numbers in the column remains the same. This is possibly due to the polarity in the distribution of the roots where this particular root is far from the rest that are in close proximity to each other. We believe that the success of our scheme in avoiding such peculiarity is in the robustness of the iterative procedure that resulted in improvement of ALL roots on equal footings whereas in other routines, it seems that the isolated root is not treated equitably and the accumulated errors are dumped into that root. Although the error is relatively small at each stage of the root-finding process in these packages, substitution in such high degree polynomial with complex coefficients does magnify the error greatly.

We believe that a rigorous assessment of the utility, accuracy and convergence properties of our routine by specialists in the field of computing software development will shed light on its utility as a viable alternative to other methods. Such assessment may lead to an improvement on the efficiency of the routine. For example, the total number of iterations, *K*, could be minimized in lieu of the iteration rounds *J* and/or the iteration steps *N*. Additionally, a better algorithm to choose the initial guessed roots may also be developed. Nonetheless, we describe here a simple but meaningful efficiency test using the polynomial example above as shown in Table IV. The Table gives a rough estimate of the convergence, accuracy and speed of the routine as explained in the Table caption. In another test that demonstrates robustness of the scheme, we used a random number generator with uniform distribution to construct the complex coefficients of a polynomial with degree 99 [7]. The real and imaginary parts of the coefficients were distributed randomly in the range $[-5,+5]$. Convergence is reached for several choices of initial guess values and number of iterations and with accuracy superior to that of Mathcad.

To illustrate how the accuracy of the individual roots develop with iteration, we plot $|P_n(\bar{z}_m)|$ in Figure 2 corresponding to $N=1$ of Table III for each iteration step until the desired accuracy is reached after 17 rounds. Figure 3 shows the same for the two cases corresponding to $N=17$ and $N=30$ in Table III but for clarity we show only the first 10 roots. Additional figures are available upon request from the authors.

Finally, we make the following relevant observations:

1. Convergence improves if we place the initial guessed roots separated on an outward or inward spiral. Table V is a reproduction of Table III with the exception that instead of placing the guessed roots on the unit circle, we fan them out on a $2\pi$-spiral starting from $x=0.5$ to $x=1.5$. It is evident that this configuration leads to more possibilities for reaching convergence with less total number of iterations. However, numerical divergences may occur for very large degrees unless we place the initial guessed roots on the unit circle.
2. Accuracy improves if after each iteration sequence ($N$ iteration steps), we reorder the resulting intermediate roots $\{\bar{z}_m\}_{m=0}^{n-1}$ in ascending order of accuracy from least accurate (largest value of $|P_n(\bar{z}_m)|$) to most accurate (smallest value of $|P_n(\bar{z}_m)|$). This is because the iteration sequence in the factorization equation (4) will result in better improvement on the value of $z_0$ if the rest $\{\hat{z}_m\}_{m=1}^{n-1}$ are the most accurate in the set.
3. The improvement in accuracy of our results over those obtained by Mathcad and Mathematica becomes more evident if the distribution of the zeros is such that one or more are located far enough from the rest in the complex plane.
4. All results obtained by our routine in this section were under the condition of "hard" convergence where all $\{|P_n(\bar{z}_m)|\}_{m=0}^{n}$ are required to be smaller than a given accuracy threshold. However, "soft" convergence allows for more possibilities to obtain accurate results if we demand that only the average of the set $\{|P_n(\bar{z}_m)|\}_{m=0}^{n}$ remains constant. On the other hand, accuracy increases in soft convergence if we include only the most accurate portion of the set $\{|P_n(\bar{z}_m)|\}_{m=0}^{n}$ in the averaging rather than the whole set.
5. For a given set of coefficients $\{a_m\}_{m=0}^{n-1}$ there exist an optimum number of iterations $N=N_c$ to achieve maximum accuracy in the values of the computed roots. Increasing $N$ beyond $N_c$ will result in diminished accuracy.

6. Our scheme does not accept initial guessed roots $\{\hat{z}_m\}_{m=0}^{n-1}$ that are degenerate (i.e., two or more guessed roots are equal). Otherwise, divergences will occur. For that reason, we placed them separated on a circle centered at the origin of the complex plane. On the other hand, the scheme is successful in finding degenerate roots. In Table VI, we give two examples illustrating the success of **Khandug** in handling a case of polynomials with doubly degenerate roots and another case with triply degenerate roots.
7. The hope is that further effort by professionals in software programming will result in enhancing the efficiency of our proposed scheme while maintaining the superior accuracy that it enjoys.

## Acknowledgements:

We are grateful to Al-Khandug Private School, Al-Madinah, Saudi Arabia for hospitality during the progress of this work. We also appreciate very much the assistance by Hasan Abdullah, Ahmed Jellal, Abdellatif Kamal, and Abdelhadi Bahaoui in the calculation using Mathematica and FORTRAN. We thank Saeed Al-Marzoug for producing the Mathematica results in the last three columns of Table II.

# Figures Caption

**Fig. 1**: The simplest version of the procedure **Khandug** that could be turned into a computational routine using any convenient programming language.

**Fig. 2** (color online): Plot of $\log_{10}|P_n(\overline{z}_m)|$ as it develops with iteration until the desired accuracy is reached for the individual roots of the polynomial example with $N=1$. Part a (Part b) is for the first (last) 10 roots.

**Fig. 3** (color online): Same as Figure 2 but for $N=17$ (part a) and $N=30$ (part b). For brevity, we show traces for the first 10 roots only.

# Tables Caption:

**Table I**: The polynomial coefficients $\{a_m\}_{m=0}^{n-1}$ used to generate the results in the Tables indicated.

**Table II**: Comparison of the accuracy of our iterative routine **Khandug** to those of Mathcad and Mathematica. We used the built-in function **polyroots** in Mathcad with the CM and LG methods. In Mathematica, results were obtained using **NRoots**, **Nsolve** and **Reduce**. The polynomial coefficients $\{a_m\}_{m=0}^{n-1}$ for $n=20$ are listed in the first column of Table I. The roots obtained by our iterative procedure are shown in the second column. The accuracy measure is given by the *sorted* set $\{|P_n(\overline{z}_m)|\}_{m=0}^{n-1}$ for the roots obtained by each of the three methods. The "0.0000" is a number less than $10^{-307}$. The individual roots in the second column are not in correspondence with the entries in the last three due to sorting.

**Table III**: Alternative convergent iteration sequences for finding the roots of the example polynomial to the desired accuracy shown in Table II and for the same set of initial guess roots, which are equally spaced on the unit circle in the complex plane. The number of rounds of the iteration procedure is $J$ whereas $K$ is the total number of iteration steps: $K = n \cdot N \cdot J$.

**Table IV**: Efficiency analysis of the accuracy and speed of convergence of our root-finding routine for the example given in Section 3. $r$ is the radius of the circle on which the initial guess roots are placed at equal separation. We vary $r$ from 0.2 to 2.2 in 10 steps until convergence is reached (exception is the $N=5$ case, where we took 50 steps). The minimum, maximum and average of $\{|P_n(\overline{z}_m)|\}_{m=0}^{n-1}$ (not including 0.000) are shown for each choice of $N$. The total number of iteration steps is $K = n \cdot N \cdot J$, where $J$ is the number of rounds of the iterative procedure to reach the desired accuracy given in Table II.

**Table V**: Reproduction of Table III with the exception that instead of placing the initial guess roots on the unit circle, we fan them out on a $2\pi$-spiral that starts from $x=0.5$ to $x=1.5$. Improved convergence is evident.

**Table VI**: An illustration that **Khandug** can handle degenerate roots successfully. We consider double degeneracy (left two columns) and triple degeneracy (right two columns) where the $n-2$ and $n-3$ polynomial coefficients are generated randomly then modified and augmented

using the degenerate root $1.0+0.5i$ to give $\{a_m\}_{m=0}^{n-1}$ for $n = 20$, which are shown in the second and third column of Table I, respectively.

**Table I**

| Tables II-V | Table VI Left | Table VI Right |
|---|---|---|
| 2+3i | -5.426-2.219i | -1.336+4.353i |
| -1 | 6.211-3.943i | 3.282-5.965i |
| -3-7i | 5.482+10.654i | -18.531+6.970i |
| 5 | -7.442-4.685i | 22.349-21.788i |
| 7+3i | -8.977-2.062i | -9.804+25.109i |
| 1+i | 10.569-3.929i | 4.616-20.920i |
| 4+2i | 0.618+2.994i | -1.602+15.328i |
| -5i | -0.290-4.231i | -5.823-15.892i |
| -7 | 3.074+5.383i | 18.622+18.195i |
| i | 2.193+1.668i | -18.746-16.666i |
| 2+8i | -12.608+6.221i | 17.098+6.774i |
| 2 | 8.509-14.117i | -12.847+10.881i |
| -7 | -1.680+3.445i | 5.197-2.908i |
| 8-2i | 1.769+9.728i | -7.739-17.275i |
| 6 | -5.230-9.811i | 9.944+10.947i |
| 5+4i | 7.841+3.444i | -4.703+9.778i |
| 2 | -4.597+1.812i | 1.909-13.533i |
| 3 | 2.443-3.728i | -3.245-0.506i |
| -1+i | -1.766+7.658i | 6.780+9.361i |
| -6i | -2.056-4.473i | -5.447-3.898i |

**Table III**

| N | 1 | 17 | 19 | 21 | 30 | 34 | 42 | 50 | 59 | 66 | 71 | 72 | 88 | 97 |
|---|---|---|---|---|---|---|---|---|---|---|---|---|---|---|
| J | 17 | 4 | 4 | 3 | 3 | 4 | 3 | 3 | 3 | 3 | 2 | 3 | 3 | 3 |
| K | 340 | 1360 | 1520 | 1260 | 1800 | 2720 | 2520 | 3000 | 3540 | 3960 | 2840 | 4320 | 5280 | 5820 |

**Table IV**

| N | r | J | K | Min ($10^{-16}$) | Max ($10^{-13}$) | Average ($10^{-14}$) |
|---|---|---|---|---|---|---|
| 1 | 1.00 | 17 | 340 | 3.6649 | 1.0049 | 1.8427 |
| 2 | 0.60 | 15 | 600 | 7.4378 | 1.6218 | 2.4595 |
| 3 | 1.00 | 4 | 240 | 7.3648 | 1.0049 | 1.9200 |
| 5 | 0.96 | 4 | 400 | 3.6649 | 1.6218 | 2.1839 |
| 10 | 0.80 | 4 | 800 | 3.6649 | 1.0049 | 1.9894 |
| 15 | 0.40 | 5 | 1500 | 3.6649 | 1.0049 | 2.1596 |
| 20 | 0.80 | 4 | 1600 | 3.6649 | 1.0049 | 2.0259 |
| 50 | 0.80 | 3 | 3000 | 3.6649 | 1.0049 | 1.9706 |
| 100 | 0.80 | 4 | 8000 | 3.6649 | 1.0049 | 2.1408 |
| 1000 | 0.80 | 3 | 60 000 | 3.6649 | 1.0049 | 2.0380 |

**Table V**

| N | 4 | 18 | 24 | 26 | 30 | 31 | 33 | 38 | 48 | 51 | 55 | 56 | 57 | 58 |
|---|---|---|---|---|---|---|---|---|---|---|---|---|---|---|
| J | 6 | 4 | 3 | 4 | 3 | 3 | 4 | 3 | 4 | 3 | 2 | 2 | 2 | 2 |
| K | 480 | 1440 | 1440 | 2080 | 1800 | 1860 | 2640 | 2280 | 3840 | 3060 | 2200 | 2240 | 2280 | 2320 |

**Table VI**

| $\{\overline{z}_m\}_{m=0}^{n-1}$ | Accuracy ($10^{-14}$) | $\{\overline{z}_m\}_{m=0}^{n-1}$ | Accuracy ($10^{-14}$) |
|---|---|---|---|
| 0.22531+3.69071i | 0.000 | 1.06018+0.07512i | 0.810 |
| 1.00285-0.01209i | 0.329 | 1.00000+0.50000i | 1.685 |
| 1.12943+0.51295i | 25.678 | 1.00000+0.50000i | 3.775 |
| 1.00000+0.50000i | 1.394 | 0.75346+0.64883i | 1.060 |
| 0.71096+0.67018i | 0.354 | 0.57299+0.77329i | 0.266 |
| 0.37696+0.93759i | 0.385 | 3.02792+3.40526i | 0.000 |
| -0.20426+1.16287i | 4.322 | 0.29746+0.89219i | 0.100 |
| -0.42241+1.03426i | 3.693 | -0.26861+0.79467i | 0.501 |
| -0.52434+0.38251i | 0.319 | -0.62921+0.77317i | 5.111 |
| -0.76117+0.41738i | 0.263 | -0.96672+0.50086i | 5.951 |
| -1.23840+0.34507i | 28.028 | -0.18275-0.38926i | 0.035 |
| -0.95149-0.24562i | 2.906 | -1.01130+0.19728i | 5.841 |
| -0.79489-0.50361i | 1.590 | -1.19670-0.50109i | 81.901 |
| 0.67188-0.46160i | 0.177 | -0.84469-0.79462i | 19.864 |
| -0.54518-0.75915i | 2.174 | -0.50381-0.80580i | 1.362 |
| -0.24322-1.25192i | 7.105 | 1.00000+0.50000i | 1.208 |
| 0.09298-0.99591i | 1.155 | 0.05341-0.95339i | 1.639 |
| 0.59499-0.98567i | 5.652 | 0.55165-1.05445i | 12.392 |
| 0.93554-0.46456i | 1.243 | 0.66244-0.73183i | 0.269 |
| 1.00000+0.50000i | 2.095 | 1.07137-0.43227i | 2.538 |

## Polynomial root-finding scheme "Khandug"

**Define**: $P(z) = z^n + a_{n-1}z^{n-1} + a_{n-2}z^{n-2} + \ldots + a_1 z + a_0$

**Input**:
1. The complex coefficients $\{a_m\}$, which could be read out from a file.
2. The number of iterations for each root: $N$
3. The maximum number of iteration rounds: $J$

**Start Khandug**
$n =$ the number of elements of the set $\{a_m\}$.
If $n = 1$ then:
    Return $\{-a_0\}$
    Exit routine
End If
If $n = 2$ then:
    Return $\left\{\tfrac{1}{2}\left(-a_1 \pm \sqrt{a_1^2 - 4a_0}\right)\right\}$
    Exit routine
End If
For $m = 0, 1, \ldots, n-1$
    $z_m = e^{i 2\pi(m/n)}$
End For $m$
For $m = 1, 2, \ldots, J$
    For $i = 0, 1, \ldots, n-1$
        $R = z_i$
        For $k = 1, 2, \ldots, N$
            $F = P(R)$
$$G = \begin{cases} \prod_{j=1}^{n-1}(R - z_j) & , i = 0 \\ \prod_{j=0}^{n-2}(R - z_j) & , i = n-1 \\ \left[\prod_{j=0}^{i-1}(R - z_j)\right]\left[\prod_{j=i+1}^{n-1}(R - z_j)\right] & , \text{otherwise} \end{cases}$$
            $R \leftarrow R - \dfrac{F}{G}$
        End For $k$
        $z_i = R$
    End For $i$
End For $m$
Return $\{z_m\}$
**End Khandug**

**Fig. 1**

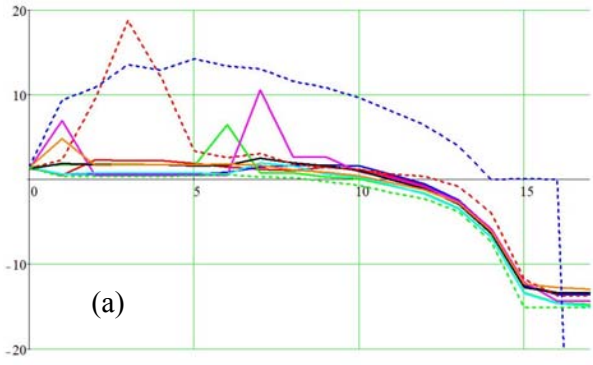
(a)

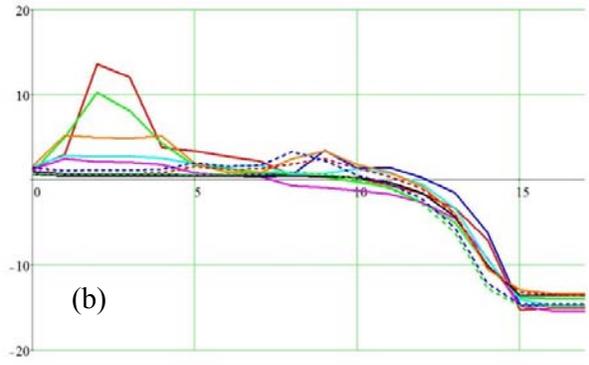
(b)

**Fig. 2**

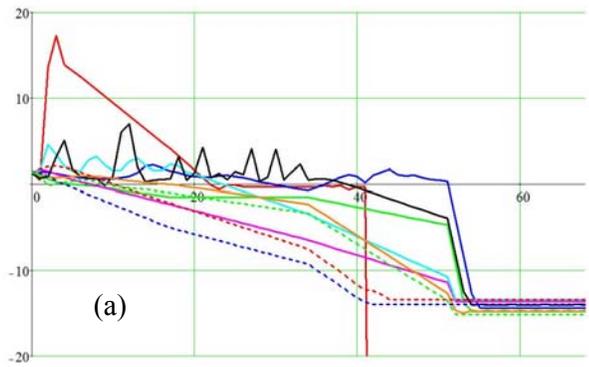
(a)

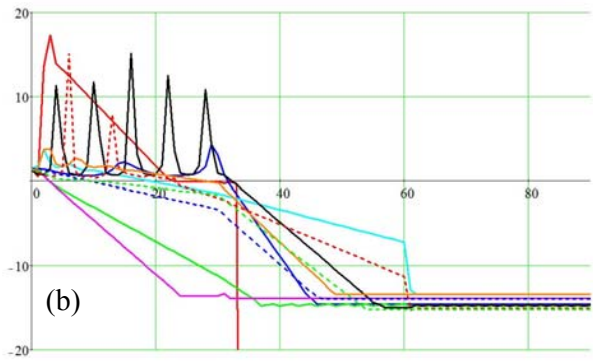
(b)

**Fig. 3**

**Table II**

| $\{|P_n(\hat{z}_m)|\}_{m=0}^{n-1}$ | $\{\bar{z}_m\}_{m=0}^{n-1}$ | Khandug ($10^{-14}$) | polyoots (CM) ($10^{-13}$) | polyroots (LG) ($10^{-5}$) | NRoots ($10^{-14}$) | NSolve ($10^{-14}$) | Reduce ($10^{-14}$) |
|---|---|---|---|---|---|---|---|
| 29.155 | -0.09074+5.81549i | 0.0000 | 0.0157 | $1.5590 \times 10^{-10}$ | 0.0099 | 0.0099 | 0.0099 |
| 22.163 | 0.89105+0.21904i | 0.0366 | 0.0673 | $2.1105 \times 10^{-9}$ | 0.0351 | 0.0570 | 0.0555 |
| 26.980 | 0.54897+0.52784i | 0.0744 | 0.0703 | 0.0291 | 0.0888 | 0.0808 | 0.0888 |
| 18.675 | 0.74874+0.81131i | 0.0828 | 0.0735 | 0.0317 | 0.0890 | 0.0993 | 0.0890 |
| 25.491 | 0.34122+1.15526i | 0.1038 | 0.0742 | 0.0322 | 0.1831 | 0.2749 | 0.2749 |
| 17.029 | -0.05401+0.89515i | 0.1559 | 0.0955 | 0.0337 | 0.2749 | 0.4441 | 0.4441 |
| 25.958 | -0.43302+1.16667i | 0.1675 | 0.1971 | 0.0358 | 0.4441 | 0.4441 | 0.7022 |
| 20.517 | 0.56340-0.21258i | 0.1691 | 0.5974 | 0.0449 | 0.4441 | 0.6779 | 0.8882 |
| 27.224 | -0.68143+0.72467i | 0.1939 | 0.6239 | 0.0464 | 0.7324 | 0.7536 | 1.0049 |
| 23.756 | -0.67024+0.23777i | 0.2726 | 0.7419 | 0.1612 | 1.0049 | 1.4648 | 1.2081 |
| 24.042 | -1.19293+0.33327i | 0.4453 | 0.9601 | 0.5237 | 1.6593 | 1.6593 | 2.1316 |
| 9.436 | -0.77633-0.08160i | 1.0893 | 1.8591 | 0.6895 | 2.1316 | 2.1316 | 2.6407 |
| 8.123 | -0.89380-0.36960i | 1.7588 | 1.9249 | 0.7331 | 3.1023 | 2.6407 | 3.5748 |
| 19.077 | -0.75081-0.62499i | 2.3110 | 4.5950 | 0.9291 | 3.1776 | 3.4307 | 5.1238 |
| 41.067 | -0.49219-0.87189i | 2.3832 | 5.4118 | 1.7956 | 5.3291 | 5.1238 | 5.3291 |
| 4.243 | 0.02587-1.00949i | 2.7716 | 12.3385 | 1.7994 | 10.1501 | 10.8939 | 6.7874 |
| 42.463 | 0.12709-1.15164i | 3.8264 | 20.7626 | 6.8793 | 11.5668 | 12.1625 | 16.0777 |
| 12.268 | 0.61401-0.77163i | 4.1431 | 39.8159 | 10.1299 | 12.1625 | 20.8977 | 26.2901 |
| 30.720 | 1.17894-0.49591i | 4.9770 | 48.3157 | 10.3193 | 20.8977 | 37.8393 | 37.8393 |
| 4.823 | 0.99621-0.29716i | 10.0486 | $7.6034 \times 10^{12}$ | $3.9375 \times 10^{10}$ | $2.5000 \times 10^{13}$ | $1.2500 \times 10^{13}$ | $2.5000 \times 10^{13}$ |